\documentclass[12pt,a4paper]{amsart}
\usepackage{amsmath,amssymb,amsthm,amscd,verbatim}
\usepackage{a4wide}
\usepackage{graphicx}
\usepackage{epsfig}
\usepackage{hyperref}

\newtheorem{thm}{Theorem}

\newtheorem*{thm*}{Theorem}

\theoremstyle{definition}

\theoremstyle{remark}

\renewcommand{\le}{\leqslant}

\def\longequation{$$\vcenter\bgroup\advance\hsize by -9em%
\noindent\ignorespaces\refstepcounter{equation}}%
\makeatletter%
\def\endlongequation{\egroup\eqno(\theequation)$$\global\@ignoretrue}
\makeatother

\begin{document}
\title{Coloring graphs with no induced subdivision of $K_4^+$}
\author{Louis Esperet} \address{Laboratoire G-SCOP (CNRS,
  Universit\'e Grenoble-Alpes), Grenoble, France}
\email{louis.esperet@grenoble-inp.fr}

\author{Nicolas Trotignon} \address{LIP (CNRS,
  ENS de Lyon), Lyon, France}
\email{nicolas.trotignon@ens-lyon.fr}

\thanks{The authors are partially supported by ANR Project STINT
  (\textsc{anr-13-bs02-0007}), and LabEx PERSYVAL-Lab
  (\textsc{anr-11-labx-0025}).}

\date{}
\sloppy

\begin{abstract}
Let $K_4^+$ be the 5-vertex graph obtained from $K_4$, the complete
graph on four vertices, by subdividing one edge precisely once (i.e. by
replacing one edge by a path on three vertices). We prove that if the
chromatic number of some graph $G$ is much larger than its clique
number, then $G$ contains a subdivision of $K_4^+$ as an induced subgraph.
\end{abstract}
\maketitle


Given a graph $H$, a \emph{subdivision} of $H$ is a graph obtained
from $H$ by replacing some edges of $H$ (possibly none) by paths. We say
that a graph $G$ contains an \emph{induced subdivision of $H$} if $G$
contains a subdivision of $H$ as an induced subgraph.

A class of graphs $\mathcal{F}$ is said to be \emph{$\chi$-bounded} if there is a
function $f$ such that for any graph $G \in \mathcal{F}$, $\chi(G)\le
f(\omega(G))$, where $\chi(G)$ and $\omega(G)$ stand for the chromatic number
and the clique number of $G$, respectively. 

\smallskip

Scott~\cite{Sco97} conjectured that for any graph $H$, the class of
graphs without induced subdivisions of $H$ is $\chi$-bounded, and
proved it when $H$ is a tree. But Scott's conjecture was
disproved in~\cite{pawlik2014}. Finding which graphs $H$ satisfy the
assumption of Scott's conjecture remains a fascinating question. It
was proved in~\cite{CELO} that every graph $H$ obtained from the
complete graph $K_4$ by subdividing at least 4 of the 6 edges once (in
such a way that the non-subdivided edges, if any, are non-incident),
is a counterexample to Scott's conjecture. On the other hand, Scott
proved that the class of graphs with no induced subdivision of $K_4$
has bounded chromatic number
(see~\cite{LMT12}). Le~\cite{Le17} proved that every graph
in this class has chromatic number at most 24.  If triangles are also
excluded, Chudnovsky et al.~\cite{chudLSSTK} proved that the chromatic number is at
most~3. 

In this paper, we extend the list of graphs known to satisfy Scott's
conjecture. Let $K_4^+$ be the 5-vertex graph obtained
from $K_4$ by subdividing one edge precisely once.

\begin{thm}\label{thm:main}
The family of graphs with no induced
subdivision of $K_4^+$ is $\chi$-bounded.
\end{thm}

We will need the following result of K\"uhn and Osthus~\cite{KO04}.

\begin{thm}[\cite{KO04}]\label{thm:ko}
For any graph $H$ and every integer $s$ there is an integer $d=d(H,s)$
such that every graph of average degree at least $d$ contains the
complete bipartite graph $K_{s,s}$ as a subgraph, or an induced subdivision of
$H$.
\end{thm}

\noindent \emph{Proof of Theorem~\ref{thm:main}.}
Let $k$ be an integer, let $d(\cdot,\cdot)$ be the function defined in
Theorem~\ref{thm:ko}, and let $R(s,t)$ be the Ramsey number of
$(s,t)$, i.e. the smallest $n$ such that every graph on $n$
vertices has a stable set of size $s$ or a clique of size $t$.

We will prove that every graph $G$ with no induced subdivision of
$K_4^+$, and with clique
number at most $k$, is $d$-colorable, with $d=\max(k,d(K_4^+,R(4,k)))$. The proof proceeds
by induction on the number of vertices of $G$ (the result being
trivial if $G$ has at most $k$ vertices). Observe that all induced
subgraphs of $G$ have clique number at most $k$ and do not contain any
induced subdivision of $K_4^+$. Therefore, by the induction, we can
assume that all induced subgraphs of $G$ are $d$-colorable. In
particular, we can assume that $G$ is connected.

Assume first that $G$ does not contain $K_{s,s}$ as a
subgraph, where $s=R(4,k)$. Then by
Theorem~\ref{thm:ko}, $G$ has average
degree less than $d$, and hence contains a vertex of degree at most
$d-1$. By the induction, $G-v$ has a $d$-coloring and this coloring
can be extended to a $d$-coloring of $G$, as desired.

We can thus assume that $G$ contains $K_{s,s}$ as a
subgraph. Since $G$ has clique number at most $k$, it follows from the
definition of $R(4,k)$ that $G$ contains $K_{4,4}$ as an induced
subgraph. Let $M$ be a set of vertices of $G$ inducing a complete
multipartite graph with at least two partite sets containing at least
4 vertices. Assume that among all such sets of vertices of $G$, $M$ is
chosen with maximum cardinality. Let $V_1,V_2,\ldots,V_t$ be the
partite sets of $M$.

Let $v$ be a vertex of $G$, and $S$ be a set of vertices not
containing $v$. The vertex $v$ is \emph{complete} to $S$ if $v$ is adjacent
to all the vertices of $S$, \emph{anticomplete} to $S$ if $v$ is not adjacent
to any of the vertices of $S$, and \emph{mixed} to $S$ otherwise. Let
$R$ be the vertices of $G$ not in $M$. We can assume that $R$ is
non-empty, since otherwise $G$ is clearly $k$-colorable and $k\le
d$. We
claim that:

\begin{longequation}\label{cl:1}
If a vertex $v$ of $R$ has at least two neighbors in some set $V_i$, then it is not
mixed to any set $V_j$ with $j\ne i$.
\end{longequation}

\medskip

Assume for the sake of contradiction that $v$ has two neighbors $a,b$
in $V_i$ and a neighbor $c$ and a non-neighbor $d$ in $V_j$, with $j\ne
i$. Then $v,a,b,c,d$ induce a copy of $K_4^+$, a contradiction. This
proves~(\ref{cl:1}).

\begin{longequation}\label{cl:2}
Each vertex $v$ of $R$ has at most one neighbor in each set $V_i$.
\end{longequation}

\medskip

Assume for the sake of contradiction that some vertex $v\in R$ has two
neighbors $a,b$ in some set $V_i$. Then by~(\ref{cl:1}), $v$ is
complete or anticomplete to each set $V_j$ with $j\ne i$. Let
$\mathcal{A}$ be the family of sets $V_j$ to which $v$ is
anticomplete, and let $\mathcal{C}$ be the family of sets $V_j$ to which $v$ is
complete. If $\mathcal{A}$ contains at least two elements, i.e. if $v$ is
anticomplete to two sets $V_j$ and $V_{j'}$ then by taking $u\in V_j$
and $u' \in V_{j'}$, we observe that $v,a,b,u,u'$ induces a copy of
$K_4^+$, a contradiction. It follows that $\mathcal{A}$ contains at
most one element. 

Next, we prove that $v$ is complete to $V_i$. Assume instead that $v$
is mixed to $V_i$. If $v$ is complete to some set $V_{\ell}$
containing at least two vertices, then we obtain a contradiction
with~(\ref{cl:1}). It follows that all the
elements of $\mathcal{C}$ are singleton. By the definition of $M$,
this implies that $\mathcal{A}$ contains exactly one set $V_j$, which has
size at least 4. Let $c$ be a non-neighbor of $v$ in $V_i$, and let
$d,d'$ be two vertices in $V_j$. Then $v,a,b,c,d,d'$ is an induced
subdivision of $K_4^+$, a contradiction. We proved that $v$ is
complete to $V_i$. Hence, every set $V_j$ is either in $\mathcal{A}$
or in $\mathcal{C}$. Since $\mathcal{A}$ contains at most one element, the
graph induced by $M\cup \{v\}$ is a complete multipartite graph, with at
least two partite sets containing at least 4 elements. This
contradicts the maximality of $M$, and concludes the proof
of~(\ref{cl:2}).

\begin{longequation}\label{cl:3}
Each connected component of $G - M$ has at most one neighbor in each set $V_i$.
\end{longequation}

\medskip

Assume for the sake of contradiction that some connected component of
$G - M$ has at least two neighbors in some set $V_i$. Then there is a
path $P$ whose endpoints $u,v$ are in $V_i$, and whose internal
vertices are in $R$. Choose $P,u,v,V_i$ such that $P$ contains the least
number of edges. Note that by~(\ref{cl:2}), $P$ contains at least 3
edges. Observe also that by the minimality of $P$, the only edges in
$G$ between $V_i$ and the internal vertices of $P$ are the first and
last edge of $P$. Let $V_j$ be a partite set of $M$ with at least 4
elements, with $j\ne i$ (this set exists, by the definition of
$M$). By~(\ref{cl:2}) and the minimality of $P$, at most two vertices of $V_j$ are
adjacent to some internal vertex of $P$. Since $V_j$ contains at least
four vertices, there exist $a,b \in V_j$ that are not adjacent to any
internal vertex of $P$. If $V_i$ has at least three elements then it
contains a vertex $w$ distinct from $u,v$. As $w$ is not adjacent to
any vertex of $P$, the vertices $w,a,b$ together with $P$ induce a
subdivision of $K_4^+$, a contradiction. If $V_i$ has at most two
elements, then there must be an integer $\ell$ distinct from $i$ and $j$
such that $V_\ell$ has at least four elements. In particular, $V_\ell$
contains a vertex $c$ that is not adjacent to any internal vertex of
$P$. As a consequence, the vertices $a,c$ together with $P$ induce a
subdivision of $K_4^+$, which is again a contradiction. This
proves~(\ref{cl:3}).

\medskip

Recall that we can assume that $R$ is non-empty. An immediate
consequence of~(\ref{cl:3}) is that the neighborhood of each connected
component of $R$ is a clique. Since $G$ is connected, it follows that
it contains a \emph{clique cutset} $K$ (a clique whose deletion
disconnects the graph). Let $C$ be a connected component of $G-K$, let
$G_1=G-C$, and let $G_2$ be the subgraph of $G$ induced by $C\cup
K$. It follows from the induction that there exist $d$-colorings of
$G_1$ and $G_2$. Furthermore, since $K$ is a clique, we can assume
that the colorings coincide on $K$. This implies that $G$ is
$d$-colorable and concludes the proof of Theorem~\ref{thm:main}.
\hfill $\Box$

\bigskip

We remark that we could have used $K_{3,3}$ instead of $K_{4,4}$ in
the proof, at the expense of a slightly more detailed analysis. The
resulting bound on the chromatic number would have been
$\max(k,d(K_4^+,R(3,k)))$ instead of $\max(k,d(K_4^+,R(4,k)))$.

\medskip

\noindent {\it Acknowledgement.} The main result of this paper was proved in January 2016
during a meeting of the ANR project STINT at Saint Bonnet de
Champsaur, France. We thank the organizers and participants for the
friendly atmosphere. We also thank Alex Scott for spotting a couple of
typos in a previous version of the draft.

\end{document}